\documentclass[12pt,oneside]{article}
\usepackage{amsmath,amssymb}
\begin{document}
%

\begin{center}
{\large\bf
On the Equivalence Problem
of Generalized  Abel ODEs under the Action of the Linear Transformations Pseudogroup}
\\[5pt]
\bf Vadim V. Shurygin, Jr.\\[5pt]
{\it E-mail:} {\tt vshjr@yandex.ru, 1vshuryg@kpfu.ru}\\[5pt]
{\it Kazan Federal University, Kazan, Russia}
\end{center}

\begin{abstract}
In the present paper we establish the necessary and sufficient conditions for
two generalized Abel differential equations to be locally equivalent
under the action of the pseudogroup of linear transformations of the form
$\{x\mapsto f(x),~ y\mapsto g(x)\cdot y + h(x)\}$. These conditions are
formulated in terms of differential invariants.
\end{abstract}



\section{Introduction}

The first kind Abel differential equation is the ODE of the form
\begin{equation}
\label{Ae}
y'=a(x)y^3+b(x)y^2+c(x)y+d(x).
\end{equation}
It was introduced by Abel in the paper~\cite{Abel}.
In what follows we assume that the functions $a$, $b$, $c$, $d$ are of class $C^\infty$ and that
 $a\ne0$.

The pseudogroup $G$ of point transformations of the form
\begin{equation}
\label{G}
x\mapsto f(x),\quad y\mapsto g(x)\cdot y + h(x), \qquad f,\,g,\,h\in
C^\infty(\mathbb{R}),
\end{equation}
preserves the class of such ODEs.
The problem of equivalence of equations (\ref{Ae}) under the action of this
pseudogroup was studied in the papers of R.~Liouville~\cite{Lio} and
P.~Appell~\cite{App}.
There are two basic relative invariants (i.e., functions $F$ in coefficients of ODE~(\ref{Ae})
and their derivatives, for which the equality $F=0$ remains invariant under the action
of $G$):
$$
s_1=a,
\quad
s_3=a'b-b'a+abc-\frac29b^3-3a^2d.
$$
Starting from $s_3$ one can construct a sequence of relative invariants
by the formula
$$
s_{2n+1} = a \frac{ds_{2n-1}}{dx} - (2n - 1)s_{2n-1}
\left(
a'
+ ac - \frac13b^2\right), \quad n\ge2.
$$
Using these, one obtains a sequence of absolute invariants
(i.e., functions that remain invariant under the action of $G$):
\begin{equation}
\label{inv3}
J_1=\frac{s_5^3}{s_3^5},
\quad
J_2=\frac{s_5s_7}{s_3^4},
\quad
J_3=\frac{s_9}{s_3^3},\dots
\end{equation}

Appell proved that this sequence can also be obtained from two basic absolute invariants
$J_1$, $J_2$, by expressing $J_2$ as a function of $J_1$ and then
differentiating the result with respect
to $J_1$.

Using the transformations~(\ref{G}) one can reduce
the Abel equation to the canonical form
$$
Y'=Y^3+R.
$$
Moreover, the transformations
\begin{equation}
\label{Kh}
\widetilde{X}=K^{-2}(X+h),~~ \widetilde{Y}=KY,\qquad K,h\in\mathbb{R},~K\ne0
\end{equation}
send one canonical form to another one.

\smallskip
{\bf Theorem 1.~\cite{App,Wone}}
{\it Two Abel equations are equivalent if and only if they have the
same canonical form modulo transformations~{\rm(\ref{Kh})}.}
\smallskip

The problem of equivalence of Abel equations with non-constant invariants was
also
considered by  E.S.~Cheb-Terrab and A.D.~Roche in~\cite{CT-R}.
Note that $R=0$ if and only if $s_3=0$ (see, e.g.,~\cite{App,Wone}).
It follows that all Abel equations satisfying $s_3=0$ are equivalent.

\section{Another approach to the equivalence of Abel ODEs}

We formulate another theorem about the equivalence of Abel equations
under the action of $G$.
Every Abel equation $\mathcal E$ is a section of the 4-dimensional bundle
$$
\pi:\mathbb{R}^5\to \mathbb{R},\qquad
(x,a,b,c,d)\to x
$$
and the pseudogroup $G$ acts on these sections.
The Lie algebra $\mathfrak g$ corresponding to $G$ consists of vector fields
$$
X=\xi(x)\dfrac{\partial}{\partial x}+
(\eta(x)\cdot y +\zeta(x))\dfrac{\partial}{\partial y}.
$$
The representation of the Lie algebra $\mathfrak g$
into the Lie algebra of vector fields on $\pi$
has the form
$$
\widehat{X}=\xi\dfrac{\partial}{\partial x}
-(2\eta+\xi')a\dfrac{\partial}{\partial a}-
(\xi' b+3\zeta a+\eta b)\dfrac{\partial}{\partial b}
+(\eta'-\xi' c-2\zeta b)\dfrac{\partial}{\partial c}+
(\zeta'-\zeta c+\eta d-\xi' d)\dfrac{\partial}{\partial d}.
$$

\smallskip
{\bf Definition.}
By an {\it (absolute) differential invariant} of order $k$ of the action of
$G$ on ${\pi}$ we understand a function
$I\in J^k({\pi})$
which is constant along the orbits of the prolonged action  of $G$.

The infinitesimal version of this definition is the equality
$$
\widehat{X}^{(k)}(I)=0
$$
for all $X\in \mathfrak g$, where $\widehat{X}^{(k)}$ denotes the $k$-th prolongation
of $\widehat{X}$ to $J^k(\pi)$.
The set of all differential invariants is an algebra.

We denote the fiber coordinates in $J^k(\pi)$ by
$a'$, $a''$, \dots.
Let
$$
\displaystyle\frac{D}{Dx}=
\frac{\partial}{\partial x}+a'\frac{\partial}{\partial a}+b'\frac{\partial}{\partial b}
+c'\frac{\partial}{\partial c}+a''\frac{\partial}{\partial a'}+b''\frac{\partial}{\partial b'}
+c''\frac{\partial}{\partial c'}+\ldots
$$
denote the total derivative operator with respect to $x$.
We say that an {\it invariant derivation}
is an operator
\begin{equation}
\label{inv-der}
\nabla = A \frac{D}{Dx},\quad
A\in C^\infty(J^\infty (\pi)),
\end{equation}
which is invariant under the prolonged action of $G$.
This is equivalent to the fact that $[\nabla, \widehat{X}]=0$
for every $X\in \mathfrak g$.
The coefficient $A$ satisfies a certain PDE system, see~\cite{Lych08}.
For every differential invariant $I$, the function $\nabla I$ also is an invariant.
Obviously, for every differential invariant $J$ the operator $J\cdot\nabla$ also is
an invariant derivation.
It follows easily that  any two invariant derivations (\ref{inv-der}) are proportional.

The first nontrivial differential invariant is $J_1$, it appears in order 2.
One can easily verify that
$$
\nabla=\frac{s_1}{s_3^{2/3}}\frac D{Dx}
$$
is an invariant derivation.
Note that invariants (\ref{inv3}) satisfy the equality
$$
J_2=\nabla (J_1^{1/3})+15J_1.
$$

The submanifold $\{s_3=0\}\subset J^1(\pi)$ is a singular orbit
for the action of $G$.
We call the point $z_k\in J^k({\pi})$ {\it regular}, if
$s_1s_3\ne0$ at this point.
In what follows  we consider only
the orbits of regular points.

Let $J$ be a differential invariant such that $DJ/Dx\ne0$ on some open interval $\Delta$.
Then for every function $F$ on $\Delta$ one has
$$
\frac {DF}{Dx}=\lambda \frac {DJ}{Dx}.
$$
The coefficient $\lambda$ is called the {\it Tresse derivative} of $F$ and is denoted
by $\lambda=DF/DJ$.
The operator $D/DJ$ is an invariant derivation (see~\cite{Lych08}).
For any invariant derivation $\nabla$ one has
$\nabla = K\cdot D/DJ$ for some invariant $K$. Substituting $J$ into this equality
we see that $K=\nabla J$. Thus,
\begin{equation}
\label{nJ}
\nabla = \nabla J \cdot \frac{D}{DJ}
\end{equation}
for every invariant $J$ and invariant derivation $\nabla$.

\smallskip
{\bf Theorem 2.}
{\it The algebra of differential invariants of the action of~$G$ on $\pi$
is generated by
the invariant $J_1$ and the invariant derivation
$\nabla$.
This algebra separates  regular orbits.}
\smallskip

{\bf Proof.}
One can easily see that the $k$-th prolongation of $\widehat{X}$
depends on the $(k+1)$-jets of functions $\xi$, $\eta$, $\zeta$.
Let $\varXi_i^k$, ${H}_i^k$, ${ Z}_i^k$ denote the components
of the decomposition
$$
\widehat{X}^{(k)}=\sum_{i=0}^{k+1} \biggl(\xi^{(i)}(x)\varXi_i^k+
\eta^{(i)}(x){H}_i^k+\zeta^{(i)}(x){Z}_i^k\biggr).
$$
The vector fields $\varXi_i^k$, ${H}_i^k$, ${ Z}_i^k$, $i=0,\dots, k+1$,
generate the completely
integrable distribution on $J^k(\pi)$ and its integral submanifolds are
exactly orbits of the action of $G$.

Let $\mathcal{O}_{k}$ be an orbit in $J^{k}(\pi)$. Its projection $\mathcal{O}_{k-1} =
\pi_{k,k-1} (\mathcal{O}_{k}) \subset  J^{k-1}(\pi)$ is
an orbit in $J^{k-1}(\pi)$.
Let $z_{k-1}\in J^{k-1}(\pi)$ be a point such that $\widehat{X}^{(k)}$ is
$\pi_{k,k-1}$-vertical over it.
Since the components
of $\widehat{X}^{(k)}$ depend  on $\xi^{(k+1)}$,
$\eta^{(k+1)}$, $\zeta^{(k+1)}$,
it follows that the bundles $\pi_{k,k-1} : \mathcal{O}_{k} \to \mathcal{O}_{k-1}$
are 3-dimensional for $k \ge 2$.
Orbits in the space of 2-jets can be found by direct integration
of 12-dimensional completely integrable distribution generating by the vector
fields $\varXi_i^k$, ${H}_i^k$, ${Z}_i^k$, $i=0,1,2, 3$.

Since the bundles $\pi_{k,k-1} : J^{k}(\pi) \to J^{k-1}(\pi)$
are 4-dimensional, it follows that for $k\ge2$
there is one differential invariant of pure order $k$ and so the dimension of
algebra of differential invariants of order $\le k$ equals $k-1$.

The invariant $J_1$
generates the space of differential invariants of pure order $2$
and separates  regular orbits in $J^2(\pi)$.
Moreover, $J_1$
is linear in second order
derivatives $a''$ and $b''$ (and does not depend on $c''$, $d''$),
and the coefficient ${s_1}/{s_3^{2/3}}$ in $\nabla$
is the  function on $J^1({\pi})$.
It follows that for $k\ge 1$ the invariant $\nabla^k J_1$
is  linear in $a^{(k+2)}$, $b^{(k+2)}$
and thus generates the space of differential invariants of pure order $k+2$
and separates  regular orbits.
$\Box$

\smallskip

Consider the space
$\mathbb{R}^{2}$ with coordinates $(j_1,j_{11})$.
For every Abel ODE  $\mathcal{E}$ we define the map
$$
\sigma_\mathcal{E}: \Delta \to\mathbb{R}^{2}
$$
by
$$
j_1= J_1^\mathcal{E}, \quad j_{11}=(\nabla J_1)^\mathcal{E},
$$
where $\Delta\subset \mathbb{R}$ is an open interval and
the superscript $\mathcal{E}$ means that the invariants are
evaluated at the coefficients of $\mathcal{E}$.

Clearly, the image
$$
\Sigma_\mathcal{E} = {\rm im}(\sigma_\mathcal{E})\subset \mathbb{R}^{2}
$$
depends only on equivalence class of $\mathcal{E}$.

\smallskip
{\bf Definition.}
We say that the equation  $\mathcal{E}$ is {\it regular}
at a point $x\in \mathbb{R}$, if\\
i) 2-jets of coefficients of $\mathcal{E}$ belong to regular orbits;\\
ii) $\sigma_\mathcal{E}(\Delta)$ is a smooth curve
in $\mathbb{R}^{2}$ for some open interval $\Delta$, containing $x$;\\
iii) one of the functions $j_1$, $j_{11}$
can be chosen as a local
coordinate on $\Sigma_\mathcal{E}$.

\smallskip
{\bf Theorem 3.}
{\it Two regular equations $\mathcal{E}$ and $\overline{\mathcal{E}}$
are locally  $G$-equivalent if and only if}
\begin{equation}
\label{Se}
\Sigma_{\mathcal{E}}=\Sigma_{\overline{\mathcal{E}}}.
\end{equation}

\smallskip
{\bf Proof.}
The necessity is obvious.

Assume that~(\ref{Se}) holds.
Let us show that $\mathcal{E}$ and $\overline{\mathcal{E}}$ are equivalent.

Without loss of generality, we may suppose that
$j_1$ is a local coordinate on
$\Sigma_\mathcal{E}$.
Let
$$
J_{11}^\mathcal{E}=j_{11}^\mathcal{E}(j_1)
$$
on  $\Sigma_\mathcal{E}$ and
$$
J_{11}^{\overline{\mathcal{E}}}=j_{11}^{\overline{\mathcal{E}}}(j_1)
$$
on $\Sigma_{{\overline{\mathcal{E}}}}$.
The condition~(\ref{Se}) means that
\begin{equation}
\label{J2}
J_{11}^\mathcal{E}=J_{11}^{\overline{\mathcal{E}}}.
\end{equation}
One can see by~(\ref{nJ}) that the invariant derivation $\nabla$ is proportional
to the Tresse derivative $D/DJ_1$ with coefficient $\nabla J_1$.
It follows from (\ref{J2}) and from Theorem 2 that
the restrictions to $\mathcal{E}$ and $\overline{\mathcal{E}}$
of differential invariants of all orders
coincide.
Since the algebra of differential invariants separates  regular orbits,
it follows that $\mathcal{E}$ and $\overline{\mathcal{E}}$ are $G$-equivalent.
$\Box$

{\bf Remark.}
Note that, by the inverse function theorem, the conditions ii)-iii) in the definition of a regular equation
are equivalent to the condition that either $J_1(x)$ or $\nabla J_1(x)$ have non-zero
derivative at a point $x$, that is either $\nabla J_1(x)\ne 0$ or $\nabla^2
J_1(x)\ne0$.

\section{Equivalence of generalized Abel ODEs}

By a generalized Abel differential equation we mean an ODE of the form
$$
y'=a_k(x)y^k+a_{k-1}(x)y^{k-1}+\ldots + a_1(x)y+a_0(x), \quad a_k(x)\ne0.
$$
We suppose all functions $a_0(x)$, \dots, $a_k(x)$ to be of class $C^\infty$.
We consider the problem of equivalence of such ODEs
under the action of the pseudogroup of linear transformations for the cases $k=4$ and $k=5$.
The computations of differential invariants were performed
using the {\tt Maple} packages {\tt DifferentialGeometry} and {\tt JetCalculus}
by I.M.\,Anderson.

\subsection{Regular case for $k=4$}

We start with the case  $k=4$, that is, the equation
\begin{equation}
\label{ga4}
y'=a(x)y^4+b(x)y^3+c(x)y^2 + d(x)y+e(x), \quad a(x)\ne0.
\end{equation}
These equations can be identified with the sections
of the 5-dimensional bundle
$$
\pi:\mathbb{R}^6\to \mathbb{R},\qquad
(x,a,b,c,d,e)\to x.
$$

The representation of $\mathfrak g$ into  the Lie algebra of vector fields on $\pi$
consists of vector fields
\begin{multline*}
\widehat{X}=\xi\dfrac{\partial}{\partial x}
-(3\eta+\xi')a\dfrac{\partial}{\partial a}-
(\xi' b+4\zeta a+2\eta b)\dfrac{\partial}{\partial b}
-(\xi' c+3\zeta b+\eta c)\dfrac{\partial}{\partial c}+\\[5pt]
+(\eta'-2\zeta c-\xi' d)\dfrac{\partial}{\partial d}
+ (\zeta'-\zeta d+\eta e-\xi' e)\dfrac{\partial}{\partial e}.
\end{multline*}

{\bf Definition.}
The function $F\in C^\infty(J^k\pi)$
is called the {\it relative differential
invariant of $k$-th order}, if for all $g\in G$ there holds
$$
g^* F=\mu(g)\cdot F ,
$$
where $\mu:G\to C^\infty(J^\infty\pi)$ is
a smooth function, called the \emph{weight function}.
The function $F\in C^\infty(J^k\pi)$ is called the
{\it absolute differential
invariant  of $k$-th order}, if $g^* F=F$  for all $g\in G$.

Again, we denote the fiber coordinates in $J^k(\pi)$ by $a'$, $a''$, etc.
There are four basic relative invariants --- two of order 0 and two of order 1:
$$
\begin{array}{c}
I_0= a,\qquad
I_1=8ac-3b^2,\\[3pt]
I_2=3(ab'-a'b)+ac^2-3abd+12a^2e,
\\[3pt]
I_3=8aa'(4ac-3b^2)+ 24a^2bb'-32a^3c'-3b^5+
64a^3cd-24a^2b^2d-32a^2bc^2+ 20ab^3c.
\end{array}
$$

{\bf Definition.}
We say that the point  $z_k\in J^k(\pi)$ is {\it regular}, if
$I_0I_1=a(8ac-3b^2)$ does not vanish
at this point.
In this subsection we consider
orbits of regular points only.

First absolute invariants appear in order 1. They are
$$
J_1=\frac{I_2I_0}{I_1^2},
\qquad
J_2=\frac{I_3}{|I_1|^{5/2}}.
$$

One can verify that the invariant derivation is
$$
\nabla=\frac{I_0^2}{|I_1|^{3/2}}\frac{D}{Dx}.
$$


\smallskip
{\bf Theorem 4.}
{\it The algebra of differential invariants of the action of
$G$ is generated by  $J_1$  and $J_2$
and the invariant derivation
$\nabla$. This algebra separates  regular orbits.
}
\smallskip

{\bf Proof.}
The proof is similar to that of Theorem 2.
We just mention the differences.

One expects four differential invariants in the order $\le2$.
They are $J_1$, $J_2$, $\nabla J_1$ and $\nabla J_2$.
These invariants
are linear in first order
and second derivatives respectively.
Thus, they generate the space of differential invariants of  order $\le 2$
and separate  regular orbits.

The bundles $\pi_{k,k-1} : J^{k}(\pi) \to J^{k-1}(\pi)$
are 5-dimensional, hence for $k\ge 3$
there are two differential invariants of pure order $k$.
The dimension of
algebra of differential invariants of order $\le k$ equals $2k$.
It follows that for $k\ge 2$ the invariants $\nabla^k J_i$, $i=1,2$,
are  linear in $a^{(k+1)}$, $b^{(k+1)}$, $c^{(k+1)}$.
The latter two invariants generate the space of differential invariants of pure order $k+1$
and separate  regular orbits.
$\Box$

\smallskip

Consider the space
$\mathbb{R}^{4}$ with coordinates $(j_1,j_2,j_{11}, j_{12})$.
For every generalized Abel ODE  $\mathcal{E}$ of the form~(\ref{ga4})  define the map
$$
\sigma_\mathcal{E}: \mathbb{R}\supset \Delta  \to\mathbb{R}^{4}
$$
by
$$
j_1= J_1^\mathcal{E}, \quad j_2= J_2^\mathcal{E},
\quad
j_{11}=(\nabla J_1)^\mathcal{E},
\quad
j_{12}=(\nabla J_2)^\mathcal{E},
$$
where the superscript $\mathcal{E}$ means that the invariants are
evaluated at the coefficients of $\mathcal{E}$.

\smallskip
{\bf Definition.}
We say that the equation  $\mathcal{E}$ is {\it regular}
at a point $x\in \mathbb{R}$, if\\
i) 2-jets of coefficients of $\mathcal{E}$ belong to regular orbits;\\
ii) $\Sigma_\mathcal{E} =\sigma_\mathcal{E}(\Delta)$ is a smooth curve
in $\mathbb{R}^{4}$ for some open interval $\Delta$, containing $x$;\\
iii) one of the functions $j_1$, $j_2$, $j_{11}$, $j_{12}$
can be chosen as a local
coordinate on $\Sigma_\mathcal{E}$.

\smallskip
The proof of the following Theorem 5 is similar to that of Theorem 3.

\smallskip
{\bf Theorem 5.}
{\it Two regular equations $\mathcal{E}$ and $\overline{\mathcal{E}}$
are locally  $G$-equivalent if and only if}
$$
\Sigma_{\mathcal{E}}=\Sigma_{\overline{\mathcal{E}}}.
$$

One can check that using the transformations (\ref{G}) the equation (\ref{ga4})
may be reduced to the canonical form
$$
Y'=Y^4+R_1Y^2+R_2, \quad R_1\ne 0.
$$
The transformations
$$
\widetilde{X}=K^{-3}(X+h),~~ \widetilde{Y}=KY
$$
map one canonical form to another one.

\subsection{Singular case for $k=4$}

Let us now consider the case when $I_1$ vanishes identically,
that is, the class of ODEs (\ref{ga4}) for which
$$
8ac=3b^2.
$$

This class consists of ODEs of the form
\begin{equation}
\label{sin4}
y'=(p(x)y+q(x))^4+r(x)y+s(x), \quad p(x)\ne 0.
\end{equation}

These equations can be identified with the sections
of the 4-dimensional bundle
$$
\pi:\mathbb{R}^5\to \mathbb{R},\qquad
(x,p,q,r,s)\to x
$$
and
the representation of $\mathfrak g$ into the Lie algebra of vector fields on $\pi$
consists of vector fields
$$
\widehat{X}=\xi\dfrac{\partial}{\partial x}
-\dfrac14(3\eta+\xi')p\dfrac{\partial}{\partial p}+
\dfrac14((\eta -\xi') q-4\zeta p)\dfrac{\partial}{\partial q}
+(\eta'-\xi' r)\dfrac{\partial}{\partial r}+
(\zeta'+(\eta -\xi') s-\zeta r)\dfrac{\partial}{\partial s}.
$$

There are three basic
relative invariants:
$$
\begin{array}{c}
L_0= p,\qquad
L_1=q'p-p'q+p^2s-pqr,
\\[5pt]
L_2=p(pq''-qp'')
+6p'(p'q-pq')
+p'p( 9qr-4ps)
-5rp^2q' -p^2qr'+p^3s'
+4p^2r(qr-ps).
\end{array}
$$
In fact, they are restrictions of
relative invariants from the regular case.
%

We restrict ourselves to the orbits of the points for which $L_1$ does not vanish.
Then there is one absolute invariant of order 2
$$
J=\frac{L_2}{L_0^2\cdot |L_1|^{7/2}}.
$$
The
invariant derivation is
$$
\nabla=\frac{|L_0|^{1/2}}{|L_1|^{3/4}}\frac{D}{Dx}.
$$

The following two theorems are proved in the same manner as above.

\smallskip
{\bf Theorem 6.}
{\it The algebra of differential invariants of the action of
$G$ is generated by  $J$
and the derivation
$\nabla$. This algebra separates  regular orbits.
}


\smallskip

{\bf Theorem 7.}
{\it Two ODEs $\mathcal{E}$ and $\overline{\mathcal{E}}$
of the form~{\rm (\ref{sin4})}
are   $G$-equivalent in the neighborhood of a point $x$ if and only if
$$
J^\mathcal{E}=J^{\overline{\mathcal{E}}}, ~~
(\nabla J)^\mathcal{E}=(\nabla J)^{\overline{\mathcal{E}}}
$$
and one of the functions $J^\mathcal{E}$ or $(\nabla J)^\mathcal{E}$ has non-zero
derivative at $x$.
}

\smallskip
Note that the canonical form of the ODE (\ref{sin4}) is
$$
Y'=Y^4+\dfrac{L_1}{L_0}\exp(-r).
$$
Thus, the ODEs (\ref{sin4}) having $L_1=0$ are all equivalent to the ODE $Y'=Y^4$.

\subsection{Regular case for $k=5$}

Now we deal  with the  equations of the form
\begin{equation}
\label{ga5}
y'=a(x)y^5+b(x)y^4+c(x)y^3 + d(x)y^2+e(x)y+f(x), \quad a(x)\ne0.
\end{equation}
These equations can be identified with the sections
of the 6-dimensional bundle
$$
\pi:\mathbb{R}^7\to \mathbb{R},\qquad
(x,a,b,c,d,e,f)\to x.
$$

There are following basic relative invariants  of the action of $G$:
$$
\begin{array}{l}
K_0=a,\qquad
K_1=5ac-2b^2,\\[5pt]
K_2=4b^3-15abc+25a^2d,\\[5pt]
K_3=50ab'-50ba'+8b^2d+5acd-50abe-3bc^2+250a^2f,\\[5pt]
K_4=2500a^2da'+1500a^2bc'-2500a^3d'-1500a^2cb'+\\[3pt]
\qquad
825a^2c^2d+6000a^2b^2f-495abc^3+1440ab^2cd-3000a^2bd^2\\[4pt]
\qquad
-288b^4d-1500a^2bce+7500a^3de-15000a^3cf+108b^3c^2.
\end{array}
$$

{\bf Definition.}
We say that the point  $z_k\in J^k(\pi)$ is {\it regular}, if
$K_0K_1=a(5ac-2b^2)$ does not vanish
at this point.
In this subsection we consider
orbits of regular points only.

The action of $G$ has the following three basic absolute invariants, one of order 0
and two of order 1:
$$
J_0=\frac{K_2^2}{K_1^3},\quad
J_1=\frac{K_3K_0^2}{|K_1|^{5/2}},
\quad
J_2=\frac{K_4K_0^2}{|K_1|^{7/2}}
$$
and the invariant derivation is
$$
\nabla=\frac{K_0^3}{K_1^2}\frac{D}{Dx}.
$$

\smallskip
{\bf Theorem 8.}
{\it The algebra of differential invariants of the action of
$G$ is generated by $J_0$, $J_1$, $J_2$
and the derivation
$\nabla$. This algebra separates  regular orbits.
}

{\bf Theorem 9.}
{\it Two ODEs $\mathcal{E}$ and $\overline{\mathcal{E}}$
of the form~{\rm (\ref{ga5})}
are   $G$-equivalent in the neighborhood of a point $x$ if and only if
$$
J_0^\mathcal{E}=J_0^{\overline{\mathcal{E}}}, ~~
J_1^\mathcal{E}=J_1^{\overline{\mathcal{E}}}, ~~
J_2^\mathcal{E}=J_2^{\overline{\mathcal{E}}}, ~~
(\nabla J_0)^\mathcal{E}=(\nabla J_0)^{\overline{\mathcal{E}}},~~
(\nabla J_1)^\mathcal{E}=(\nabla J_1)^{\overline{\mathcal{E}}},~~
(\nabla J_2)^\mathcal{E}=(\nabla J_2)^{\overline{\mathcal{E}}}
$$
and at least one of the above six functions has non-zero
derivative at $x$.
}

The canonical form of the ODE (\ref{ga5}) with respect to the action of $G$ is
$$
Y'=Y^5+R_1Y^3+R_2Y^2+R_3, \quad R_1\ne 0.
$$
The transformations
$$
\widetilde{X}=K^{-4}(X+h),~~ \widetilde{Y}=KY
$$
permute the canonical forms.

\subsection{Singular cases for $k=5$}

The first singular case arises when $K_1$ vanishes, that is, when
$$
5ac=2b^2.
$$

It corresponds to the ODEs of the form
\begin{equation}
\label{sin5-1}
y'=(p(x)y+q(x))^5+r(x)y^2+s(x)y+t(x), \quad p,r\ne0.
\end{equation}

The action of $G$ on the class of ODEs~(\ref{sin5-1})
has the following basic relative invariants
$$
L_0= p,\quad
L_1=r, \quad L_2=-p'q+q'p+q^2r-pqs+tp^2,
\quad
L_3=5p'r+3prs-6qr^2-pr'.
$$
They provide two
absolute invariants of order 1:
$$
J_0=\frac{L_2\cdot p^{4/3}}{r^{5/3}},
\qquad
J_1=\frac{L_3\cdot p^{2/3}}{r^{7/3}}.
$$
The invariant derivation is
$$
\nabla=\frac{p^{5/3}}{r^{4/3}}\frac{D}{Dx}.
$$

{\bf Theorem 10.}
{\it Two ODEs $\mathcal{E}$ and $\overline{\mathcal{E}}$
of the form~{\rm (\ref{sin5-1})}
are   $G$-equivalent in the neighborhood of a point $x$ if and only if
$$
J_0^\mathcal{E}=J_0^{\overline{\mathcal{E}}}, ~~
J_1^\mathcal{E}=J_1^{\overline{\mathcal{E}}}, ~~
(\nabla J_0)^\mathcal{E}=(\nabla J_0)^{\overline{\mathcal{E}}},~~
(\nabla J_1)^\mathcal{E}=(\nabla J_1)^{\overline{\mathcal{E}}}
$$
and at least one of the above functions has non-zero
derivative at $x$.
}

The canonical form for such ODEs is
$$
Y'=Y^5+R_2Y^2+R_3.
$$

The second singular case
arises when both $K_1$ and $K_2$ vanish
and  corresponds to the ODEs of the form
\begin{equation}
\label{sin5-2}
y'=(p(x)y+q(x))^5+s(x)y+t(x), \quad p\ne0.
\end{equation}

The relative invariants are
$$
\begin{array}{c}
M_0= p,\quad
M_2=-p'q+q'p-pqs+tp^2,
\\[5pt]
M_4=p(pq''-qp'') +7p'(p'q-pq')+p^3t'
-p^2qs'-\qquad\qquad\qquad\\[5pt]
\qquad\qquad\qquad
-6p^2sq'+p'p(11qs-5pt)+5p^2s(qs-pt).
\end{array}
$$

For the subclass of ODEs for which $M_2\ne0$
there is one basic absolute invariant
$$
J=\frac{M_4}{p^{2/5}\cdot M_2^{9/5}}
$$
and the invariant derivation
$$
\nabla=\frac{p^{3/5}}{M_2^{4/5}}\frac{D}{Dx}.
$$

The canonical form for such ODEs is
$$
Y'=Y^5+R_3
$$
and the case $M_2=0$ corresponds to $R_3=0$.

\smallskip
{\bf Theorem 11.}
{\it Two ODEs $\mathcal{E}$ and $\overline{\mathcal{E}}$
of the form~{\rm (\ref{sin5-2})}
are   $G$-equivalent in the neighborhood of a point $x$ if and only if
$$
J^\mathcal{E}=J_0^{\overline{\mathcal{E}}}, ~~
(\nabla J)^\mathcal{E}=(\nabla J_0)^{\overline{\mathcal{E}}}
$$
and at least one of the above functions has non-zero
derivative at $x$.
}


\vskip1cm

\begin{thebibliography}{30}

\bibitem{Abel}
N.~Abel, {\it Pr\'ecis d'une th\'eorie des fonctions elliptiques}, J. Reine. Angew. Math. 4
(1829), 309--348.

\bibitem{App}
P.~Appell, {\it Sur les invariants de quelques \'equations diff\'erentielles},
Journal de
Math\'ematique 5, 361--423 (1889).

\bibitem{CT-R}
E.S.~Cheb-Terrab, A.D.~Roche,
{\it Abel ODEs: Equivalence and Integrable Classes},
Computer Physics Communications. 01/2000; DOI: 10.1016/S0010-4655(00)00042-4,
{\tt arXiv:math-ph/0001037}, 2000.


\bibitem{Lio}
R.~Liouville, {\it Sur une classe d'\'equations diff\'erentielles du premier ordre et les
formations invariantes qui s'y rapportent},
Comptes Rendus Acad. Sci. 105 (1887), 460--463.

\bibitem{Lych08}
V.~Lychagin, {\it Feedback Equivalence of 1-dimensional Control
Systems of the 1-st Order}, ``Geometry, topology and there
applications'', Proceedings of the Institute of mathematics of NAS of
Ukraine, 2009, 6(2), pp. 288--302.

\bibitem{Wone}
O.~Wone, {\it Geometry of the Abel Equation of the first kind},
{\tt arXiv:1401.2375 [math.DG]}, 2014.

\end{thebibliography}
\end{document}